\magnification = 1200
\input amssym.def
\input amssym.tex
\def \qed {\hfill $\square$}
\def \R {\Bbb R}
\def \Z {\Bbb Z}
\def \ms {\medskip}
\def \msi {\medskip\noindent}
\def \ssi {\smallskip\noindent}

\def \bsi {\bigskip\noindent}
\def \pari {\par\noindent}

\def \sm {\setminus }

\def \wt {\widetilde }

\overfullrule=0pt
\def \dist{\mathop{\rm dist}\nolimits}
\def \diam{\mathop{\rm diam}\nolimits}

\centerline{A CONSTRUCTIVE PROOF OF THE ASSOUAD EMBEDDING THEOREM}
\par
\centerline{WITH BOUNDS ON THE DIMENSION}

\bigskip
\centerline { Guy DAVID and Marie SNIPES}

\bigskip\noindent
{\bf R\'{e}sum\'{e}.}
On red\'{e}montre de mani\`{e}re constructive un r\'{e}sultat de 
Naor and Neiman (\`{a} para\^{\i}tre,  Revista Matematica 
Iberoamercana),
qui dit que pour tout espace m\'{e}trique doublant $(E,d)$, il existe
$N \geq 0$, qui ne d\'{e}pend que de la constante de doublement, tel que
pour tout exposant $\alpha \in ]1/2,1[$, il existe une application
bilipschitzienne $F$ de $(E,d^{\alpha})$ \`{a} valeurs dans $\R^N$.

\medskip \noindent
{\bf Abstract.}
We give a constructive proof of a theorem of Naor and Neiman,
(to appear, Revista Matematica Iberoamercana), 
which asserts that if $(E,d)$ is a doubling metric space,
there is an integer $N > 0$, that depends only on the metric doubling
constant, such that for each exponent $\alpha \in (1/2,1)$, we can find
a bilipschitz mapping $F = (E,d^{\alpha}) \to \R^N$.


\bigskip
The purpose of this paper is to give a simpler and constructive
proof of a theorem proved by
Naor and Neiman [NN], 
which asserts that if $(E,d)$ is a doubling metric space,
there is an integer $N > 0$, that depends only on the metric doubling
constant, such that for each exponent $\alpha \in (1/2,1)$, we can find
a bilipschitz mapping $F = (E,d^{\alpha}) \to \R^N$.

Here $\R^N$ is equipped with its Euclidean metric,
the snowflake distance $d^\alpha$ is simply defined by
$d^\alpha(x,y) = d(x,y)^\alpha$ for $x, y \in E$,
and metrically doubling means that there is 
an integer $C_0 \geq 1$
such that for every $r > 0$, every (closed) ball of radius $2r$
in $E$ can be covered with no more than $C_0$ balls of radius $r$.
We call $C_0$ a metric doubling constant for $(E,d)$.

\ssi{\bf Remark 1.}
Notice that in a doubling metric space with doubling constant $C_0$, 
every ball of radius $2^m r$ can be covered
with $C_0^m$ balls of radius $r$. More generally, for $\lambda > 0$,
every ball of radius $\lambda r$ can be covered by
$C_0 \lambda^{N_0}$ balls of radius $r$, where $N_0=\log_2 C_0$.
To see this replace $\lambda$ by the next power of $2$.

A consequence of Naor and Neiman's main result in [NN] 
is the following.

\ms\proclaim Theorem 2.
For each $C_0 \geq 1$, there is an integer $N$ and,
for $1/2 < \alpha < 1$, a constant $C = C(C_0,\alpha)$
such that if $(E,d)$ is a metric space that admits the
metric doubling constant $C_0$, we can find an injection
$F : E \to \R^N$ such that
$$
C^{-1} d(x,y)^\alpha \leq |F(x)-F(y)| \leq C d(x,y)^{\alpha}
\leqno (3)
$$
for $x, y \in E$.

\ms
The theorem of Naor and Neiman is more general,
but its proof is also more complicated.
If we did not ask $N$ to be independent of $\alpha$,
this would just be the usual Assouad embedding theorem
from [A]. Here we do not try to get optimal values of
$N$ and $C$; see [NN] for better control of these constants, and for
more information on the context. There is nothing special
about the constant $1/2$ in our statement; we just
do not want to consider the case when $\alpha$
is close to $0$, for which the dimension independence fails.
Indeed, even when $E= \R$, with the usual distance,
one needs many dimensions to construct an $\alpha$-snowflake
with $\alpha$ small.


The proof will be a rather simple modification of the
standard proof of Assouad's theorem, but it took
a surprising amount of energy to the authors, plus the
knowledge of the fact that the result is true, to
make it work.

Rather than using a probabilistic proof, we use an
adaptive argument, and work at small relative scales to
use the fact that there is a lot of space in $\R^N$;
the difficult part (for the authors) was to realize that
using a very sparse collection of scales would not kill
the argument, and instead helps control the residual
terms. The constant $C(C_0,\alpha)$ gets very large (when $\alpha$
gets close to $1$), but this is expected.

The issue of giving a simpler proof of Theorem 2
came out in the AIM workshop on Mapping Theory in
Metric Spaces in Palo Alto (January 2012);
the authors wish to thank Mario Bonk for asking the question,
other participants of the task force, and in particular
U. Lang, C. Smart, J. Tyson, 
for very interesting discussions and their patience
with at least one wrong proof, and L Capogna, J. Tyson,
and S. Wenger who organized the event.

\ms
We need some notation before we start the proof.
Here the doubling metric space $(E,d)$ is fixed,
and $B(x,r)$ will denote the closed ball with center $x$
and radius $r$. We may assume that $\alpha$ is close to $1$,
because otherwise we can use the standard Assouad result and proof.

We shall use a small parameter $\tau > 0$, with
$\tau \leq 1-\alpha$, and work at the scales
$$
r_k = \tau^{2k}, k \in \Z.
\leqno (4)
$$
We first prove Theorem~2 in the case that $E$ has finite diameter;
 this allows us to start at
$k = k_0$, where $k_0$ is such that $r_{k_0} \geq \diam(E)$.

For each $k \geq k_0$, select a maximal collection $\{ x_j \}$,
$j\in J_k$, of points of $E$, with $d(x_i,x_j) \geq r_k$ for $i \neq j$.
Thus, by maximality
$$
E \i \bigcup_{j\in J_k} B(x_j,r_k),
\leqno (5)
$$
For convenience, we write $B_j = B(x_j,r_k)$ and
$\lambda B_j = B(x_j,\lambda r_k)$ for $\lambda > 0$.

Letting $N(x)$ denote the number of indices $j\in J_k$ such that
$d(x_j,x) \leq 10r_k$, we now check that
$$
N(x) \leq C_0^5 \ \hbox{ for } x\in E.
\leqno (6)
$$
Indeed, we can cover $B(x,10r_k)$ with fewer than $C_0^5$
balls $D_l$ of radius $r_k/3$. Each $D_l$ contains
at most one $x_j$ (because $d(x_i,x_j) \geq r_k$ for
$i \neq j$). Because all the $x_j$ that lie in $B(x,10r_k)$
are contained in some $D_l$, (6) follows.

From inequality (6) and the assumption (for the present case) that $E$ is bounded,
we conclude $J_k$ is finite. Note, however, that even if $E$ were unbounded,
we could use the preceding discussion to construct the collection
$\{ x_j \}$ without using the axiom of choice.

Set
$\Xi = \big\{ 1, 2, \ldots C_0^5 \big\}$ (a set of colors).
Now enumerate $J_k$, and for $j\in J_k$ let $\xi(j)$ be
the first color not taken by an earlier close neighbor; specifically, choose $\xi(j)$ to be the first 
color not used by an earlier $i\in J_k$ such that $d(x_i,x_j) \leq 10r_k$.  By construction we have
$$
\xi(i) \neq \xi(j) \ \
\hbox{ for $i, j \in J_k$ with $i \neq j$ and }
d(x_i,x_j) \leq 10r_k.
\leqno (7)
$$

Finally, for each $\xi \in \Xi$, define the set $J_k(\xi) := \big\{ j \in J_k \, ; \, \xi(j)=\xi \big\}$. 
Thus
$$
d(x_i,x_j) > 10r_k
\ \hbox{ for } i, j \in J_k(\xi) \hbox{ such that } i \neq j.
\leqno (8)
$$

\ms
For each $j\in J_k$, set
$\varphi_j(x) = \max\{0, 1 - r_k^{-1}\dist(x,B_j)\}$.
The formula is not important; we just want to make sure that
$$
0 \leq \varphi_j(x) \leq 1 \ \hbox{ everywhere,}
\leqno (9)
$$
$$
\varphi_j(x) = 1 \ \hbox{ for } x\in B_j,
\leqno (10)
$$
$$
\varphi_j(x) = 0 \ \hbox{ for } x\in E \sm 2B_j,
\leqno (11)
$$
and
$$
\varphi_j \hbox{ is Lipschitz, with }
||\varphi_j||_{lip} \leq r_k^{-1}.
\leqno (12)
$$

\ms
We continue with the non-surprising part of the construction.
For each $\xi\in \Xi$, we will construct two mappings:
$F^\xi : E \to \R^M$ and a slightly modified version $\wt F^\xi : E \to \R^M$, where $M$ is a very large integer
depending only on the metric doubling constant.
Our final mapping $F:E\to\R^N$,
will be the direct product of these $2C_0^5$ mappings. Thus, the dimension $N$ is $2C_0^5 M$, which
can probably be improved.

We decide that $F^\xi$ will be of the form:
$$
F^\xi(x) = \sum_{k \geq k_0} r_k^\alpha f_k^\xi(x),
\leqno (13)
$$
where
$$
f_k^\xi(x) = \sum_{j\in J_k(\xi)} v_{j} \varphi_j(x),
\leqno (14)
$$
and with vectors $v_j \in \R^M$ that will be carefully chosen later.  $\wt F$ will take on the same form as $F$, but with a different choice of the vectors $\{v_j\}$.  However, for both functions we will choose the $v_j$ inductively, and so that
$$
v_j \in B(0,\tau^2) \i \R^M,
\leqno (15)
$$
with the same very small $\tau > 0$ as in the definition
of $r_k = \tau^{2k}$ above; $\tau$ will be chosen near the end.
With this choice, we immediately see that
$$
|| f_k^\xi ||_\infty \leq \tau^2
\leqno (16)
$$
because the $\varphi_j$, $j\in J_k(\xi)$,
have disjoint supports by (8) and (11);
hence the series in (14) converges.  Moreover, if we set
$$
F_k^\xi(x) = \sum_{k_0 \leq \ell \leq k} r_\ell^\alpha f_\ell^\xi(x),
\leqno (17)
$$
we get that
$$
|| F^\xi-F_k^\xi ||_\infty
\leq \sum_{\ell > k} r_\ell^\alpha \tau^2
= r_{k+1}^\alpha \tau^2 \sum_{\ell \geq 0} \tau^{2\ell\alpha}
\leq 2\tau^2 r_{k+1}^\alpha
\leqno (18)
$$
(because $r_l = \tau^{2l}$ and $\tau^{2\alpha} < 1/2$
when $\tau$ is small).
Also, the Lipschitz norm of $f_k^\xi$ is
$$
|| f_k^\xi ||_{lip} \leq \tau^2 r_k^{-1}
\leqno (19)
$$
by (12) and because the $\varphi_j$ are supported in
disjoint balls; we sum brutally and get that
$$\eqalign{
|| F_k^\xi ||_{lip}
&\leq \sum_{\ell \leq k} r_\ell^\alpha || f_\ell^\xi ||_{lip}
\leq \tau^2 \sum_{\ell \leq k} r_\ell^{\alpha-1}
= \tau^2 r_k^{\alpha-1} \sum_{\ell \leq k} \tau^{2 (\ell-k )(\alpha-1)}
\cr&
= \tau^2 r_k^{\alpha-1} (1-\tau^{2(1-\alpha)})^{-1}
}\leqno (20)
$$
by (17). We take $\tau \leq 1-\alpha$
(many other choices would do, the main point is to
have a control in (22) below by a power of $\tau$, which
could even be negative); then
$$
\ln(\tau^{2(1-\alpha)}) = 2(1-\alpha) \ln(\tau)
= -2 (1-\alpha) \ln\left({1 \over\tau}\right)
\leq -2 \tau \ln\left({1 \over\tau}\right);
\leqno (21)
$$
we exponentiate and get that
$\tau^{2(\alpha-1)}
\leq e^{-2 \tau \ln({1 \over\tau})}
\leq 1 -\tau \ln({1 \over\tau})$
if $\tau$ is small enough, hence
$1-\tau^{2(\alpha-1)} \geq \tau \ln({1 \over\tau})$,
and finally
$$
|| F_k^\xi ||_{lip} \leq {\tau \over \ln({1 \over\tau})}
\, r_k^{\alpha-1} \leq r_k^{\alpha-1}
\leqno (22)
$$
(again if $\tau$ is small enough; for example, $\tau<1/2$ works).

\ms
We now describe how to choose the
vectors $v_j$, $j\in J_k$, so that the differences
$|F_k^\xi(x)-F_k^\xi(y)|$ will be as large as possible (toward proving that $(F_k^\xi)^{-1}$ is bilipschitz).

Fix $k \geq k_0$,
suppose that the $F_{k-1}^\xi$ were
already constructed, and fix $\xi \in \Xi$.
Put any order $<$ on the finite set $J_k(\xi)$.
We shall construct $F_{k}^\xi$ with the order $<$ and $\wt F_{k}^\xi$ with the reverse order.

Choose the $v_j$ for $F_{k}^\xi$ according to the order $<$.
Recall that we defined
$$
F_{k}^\xi(y) = F_{k-1}^\xi(y) + r_k^\alpha f_k^\xi(y)
= F_{k-1}^\xi(y) + r_k^\alpha \sum_{i\in J_k(\xi)}  v_{i} \varphi_i(y)
\leqno (23)
$$
in (13) and (14); for each $j\in J_k(\xi)$, we shall also consider
the partial sum $G_{k,j}^\xi$ defined by
$$
G_{k,j}^\xi(y) = F_{k-1}^\xi(y) + r_k^\alpha
\sum_{i\in J_k(\xi) \, ; \, i<j}  v_{i} \varphi_i(y),
\leqno (24)
$$
which we therefore assume to be known when we choose $v_j$.

\ms\proclaim Lemma 25.
For each $j \in J_k(\xi)$, we can choose $v_j \in B(0,\tau^2)$
so that
$$
|F_{k}^\xi(x) - G_{k,j}^\xi(y)| \geq \tau^3 r_k^\alpha
\ \  \hbox{ for $x\in B_j$ and $y\in B(x_j,10\tau^{-2}r_k) \sm 2B_j$.}
\leqno (26)
$$

\ms
The extra room in $\R^M$
will be used to give lots of different choices of $v_j$.
Observe that for $x\in B_j$, $\varphi_j(x) = 1$ and the
other $\varphi_i(x)$ are all null (because $\varphi_i$
is supported in $2B_i$ by (11), and $2B_i$ never meets $B_j$
by (8)). Thus
$$
F_{k}^\xi(x) = F_{k-1}^\xi(x) + r_k^\alpha f_k^\xi(x)
= F_{k-1}^\xi(x) + r_k^\alpha v_j
\leqno (27)
$$
by (13) and (14). By (22),
$|| F_k^\xi ||_{lip} \leq r_k^{\alpha-1}$,
but the proof of (22) also yields
$$
|| G_{k,j}^\xi ||_{lip} \leq r_k^{\alpha-1}
\leqno (28)
$$
(we just add fewer terms).
We shall use this to replace $B_j$
and $B(x_j,10\tau^{-2}r_k) \sm 2B_j$ with discrete
sets.
Set $\eta = \tau^3 r_k$, and pick an
$\eta$-dense set $X$ in $B_j$, and
an $\eta$-dense set $Y$ in $B(x_j,10\tau^{-2}r_k) \sm 2B_j$.
We shall soon prove that we can choose $v_j$ so that
$$
|F_{k}^\xi(x') - G_{k,j}^\xi(y')| \geq 3\tau^3 r_k^\alpha
\ \hbox{ for $x'\in X$ and } y'\in Y,
\leqno (29)
$$
and let us first check that the lemma will follow.

Notice that for $x\in B_j$, we can find
$x'\in X$ such that
$$
|F_{k}^\xi(x')-F_{k}^\xi(x)|
\leq || F_{k}^\xi ||_{lip} \, \eta
\leq r_k^{\alpha-1} \cdot \tau^3 r_k
= \tau^3 r_k^\alpha,
\leqno (30)
$$
and similarly, for $y \in B(x_j,10\tau^{-2}r_k) \sm 2B_j$
we can find $y'\in Y$ such that
$$
|G_{k,j}^\xi(y)-G_{k,j}^\xi(y')|
\leq || G_{k,j}^\xi ||_{lip} \, \eta
\leq \tau^3 r_k^\alpha.
\leqno (31)
$$
Then (26) for $x$ and $y$ follows from (29), as needed.

So we want to arrange (29).
First we bound $|X|$, the number of elements in $X$.
By Remark 1, we can cover $B_j = B(x_j,r_k)$ by
$C_0 (2\tau^{-3})^{N_0}$ balls of radius $\eta/2$;
we just keep those that meet $B_j$, pick an element
of $B_j$ in each such ball, and get an $\eta$-dense net
$X$, with $|X| \leq C_0 (2\tau^{-3})^{N_0}$.
Similarly, we can find $Y$ so that
$$
|Y| \leq C_0 \left(2 \, {10 \tau^{-2} r_k \over \eta}\right)^{N_0}
= C_0 \left({20 \tau^{-2} r_k \over \tau^3 r_k}\right)^{N_0}
= C_0 (20\tau^{-5})^{N_0}.
\leqno (32)
$$
The total number of pairs $(x',y')$ for which we have to check (29)
is thus
$$
|X||Y| \leq C_0^2 (40 \tau^{-8})^{N_0};
\leqno (33)
$$
worse estimates on $|| F_k^\xi ||_{lip}$ and $|| G_{k,j}^\xi ||_{lip}$
above would have yielded worse powers of $\tau$, but we would not care.

Now pick a maximal finite set $V$ in $B(0,\tau^2) \i \R^M$,
whose points lie at mutual distances at least $7 \tau^3$
from each other. For each pair $(x',y')$ as above,
the different choices of $v_j \in V$ yield the same
value of $G_{k,j}^\xi(y')$ (because $G_{k,j}^\xi$
does not depend on $v_j$, by (24)), and values of
$F_{k}^\xi(x')$ that differ by at least $7 \tau^3 r_k^\alpha$,
by (27). Thus (29) for this pair $(x',y')$ cannot fail for
more than one choice of $v_j \in V$, and it is now enough
to show that $V$ has more than $|X||Y|$ elements.
Taking $C_M$ to be the doubling constant of $\R^M$, it follows from Remark~1 that 
$|V| \geq C_M (1/7\tau)^M$,
which is indeed
larger than $|X||Y|$ if $M > 8 N_0$ and
$\tau$ is small enough (depending on $M$).
This completes our verification of (29); as noted earlier, Lemma~25
follows.
\qed

\ms
For each color $\xi$, choose the vectors $v_j$, and hence the
mappings $F_k^{\xi}$, as in Lemma 25. Also
define a second version $\wt F_k^{\xi}$ of $F_k^{\xi}$
using the opposite order on $J_k(\xi)$.
Define $F_k : E \to \R^N$ (with $N = 2C_0^5 M$) to be the direct product of the $2C_0^5$ maps
$F_k^{\xi}$ and $\wt F_k^{\xi}$ for all colors $\xi\in \Xi$.
Finally
set $F = \lim_{k \to +\infty} F_k$. We are ready
to check that $F$ is bilipschitz.

\ms\proclaim Lemma 34. We have that
$$
{\tau^5 \over 8} \, d(x,y)^\alpha
\leq |F(x)-F(y)| \leq 5N \tau^{-2(1-\alpha)} d(x,y)^\alpha
\ \hbox{ for } x, y \in E.
\leqno (35)
$$

\ms
Let $x, y \in E$ be given; we may assume that
$x \neq y$. Let $k$ be such that
$$
4r_k \leq d(x,y) \leq 4r_{k-1} = 4 \tau^{-2} r_k,
\leqno (36)
$$
where the last part comes from (4).
Then $r_k \leq 4r_k \leq d(x,y) \leq \diam(E)
\leq r_{k_0}$ by our definition of $k_0$,
and so $k \geq k_0$.
By (22) and (36),
$$\eqalign{
|F_k^\xi(x)-F_k^\xi(y)|
&\leq || F_k^\xi ||_{lip} \, d(x,y)
\leq r_k^{\alpha-1} d(x,y)
\cr&
\leq \Big({d(x,y) \over 4\tau^{-2}}\Big)^{\alpha-1} d(x,y)
= (4\tau^{-2})^{1-\alpha} d(x,y)^\alpha.
}\leqno (37)
$$
But (18) also says that
$$\eqalign{
\big||F^\xi(x)-F^\xi(y)| - |F_k^\xi(x)-F_k^\xi(y)|\big|
&\leq 2 || F^\xi-F_k^\xi ||_\infty
\leq 4\tau^2 r_{k+1}^\alpha
= 4 \tau^2 \tau^{2\alpha} r_k^\alpha
\cr&
\leq 2\tau^2 d(x,y)^\alpha
}\leqno (38)
$$
and (37) and (38) (and similar estimates for
the $\wt F_k^{\xi}$) give the upper bound in (35).

For the lower bound, notice that by (5) we can find
$j\in J_k$ such that $x\in B_j$.
We shall just consider the color $\xi \in \Sigma$ such that
$j\in J_k(\xi)$, and distinguish between two cases.

We'll need to know that
$$
y\in B(x_j,10\tau^{-2}r_k) \sm 2B_j.
\leqno (39)
$$
That $y\in B(x_j,10\tau^{-2}r_k)$ follows from (36), because
$d(x,x_j) \leq r_k$ since $x\in B_j$. Moreover, if $y\in 2B_j$,
then $d(x,y) \leq d(x,x_j) + d(x_j,y) \leq 3 r_k$,
which would contradict (36). So (39) holds.

First assume that $y\in 2B_i$ for some $i\in J_k(\xi)$.
Then $i \neq j$, by (39).
Let us assume that $i<j$; otherwise, we would use
$\wt F_k^{\xi}$ instead of $F_k^{\xi}$ in the following calculations.
Recall that all the $\varphi_l(y)$, $l \neq i$,
are equal to $0$, by (11) and (8). Then (23) and (24)
yield
$$
F_{k}^\xi(y)
= F_{k-1}^\xi(y) + r_k^\alpha  v_{i} \varphi_i(y)
= G_{k,j}^\xi(y)
\leqno (40)
$$
(with $F_{k-1}^\xi(y) = 0$ if $k=k_0$).
By (39), we can apply (26), which says that
$$
|F_{k}^\xi(x) - F_{k}^\xi(y)|
= |F_{k}^\xi(x) - G_{k,j}^\xi(y)|
\geq \tau^3 r_k^\alpha.
\leqno (41)
$$
We then combine this with (38) and get that
$$\eqalign{
|F^\xi(x) - F^\xi(y)|
&\geq |F_{k}^\xi(x) - F_{k}^\xi(y)|- 4 \tau^2 \tau^{2\alpha} r_k^\alpha
\geq \tau^3 r_k^\alpha - 4\tau^2 r_{k+1}^\alpha
\cr&
= \tau^3 r_k^\alpha (1-4\tau^{-1} \tau^{2\alpha})
\geq {\tau^3 r_k^\alpha \over 2}
}\leqno (42)
$$
by (4) and because we can take $\alpha > 2/3$ and
$\tau$ small. Now
$$
|F(x)-F(y)| \geq |F^\xi(x)-F^\xi(y)|
\geq {\tau^3 r_k^\alpha \over 2}
\geq {\tau^3 \over 2}\, \Big({d(x,y) \over 4\tau^{-2}}\Big)^\alpha
\geq {\tau^5 \over 8} \, d(x,y)^\alpha
\leqno (43)
$$
by (36). This proves (35) when $y\in 2B_i$ for some $i\in J_k(\xi)$.
If not, all the $\varphi_i(x)$ vanish by (11), and
$F^\xi_k(y) = G^\xi_{k,j}(y) = F^\xi_{k-1}(y)$ for all $j$,
by (23) and (24). That is, (40) still holds (with $j$ again chosen so that $x\in B_j$), 
and we can continue just as in the previous case.
Lemma 34 follows.
\qed

\ms
This completes the proof of Theorem 2 for bounded $E$.

  Now suppose $E$ is an unbounded metric space with doubling constant $C_0$.
Fix an origin $x_0$, and apply the construction
above to the sets $E_m = E \cap B(x_0,2^m)$.

The set $E_m$ is itself doubling, the doubling constant
$C_0^2$.  To see this note that if $x\in E_m$ and $r>0$, we can cover
$E_m \cap B(x,2r)$ with $C_0^2$ balls of radius $r/2$,
which (when they meet $E_m$) we can replace with balls of radius $r$ whose centers are in $E_m$.

We get from the proof above a mapping $F_m$ such that
$$
C^{-1} d(x,y)^\alpha \leq |F_m(x)-F_m(y)| \leq C d(x,y)^{\alpha}
\leqno (44)
$$
for $x, y \in E_m$, where $C$ depends on
$C_0$ and $\alpha$ but not on $m$.
We may assume that $F_m(x_0) = 0$, after possibly adding a constant, 
which would not destroy inequality (44).

Now define for each $k \in \Bbb Z$, a maximal collection $\{ x_j \}\subset E$, $j\in J_k$, 
with $d(x_i, x_j)\geq r_k$ for $i\neq j$.
Although $E$ is unbounded, each $J_k$ is still at most countable.

Notice also that for each $x_j$, the sequence $\{ F_m(x_j) \}$
is bounded (by (44) and because $F_m(x_0) = 0$); hence we can
extract a subsequence $\{ m_j \}$, so that the sequence $F_{m_j}(x_j)$
converges for each $x_j$ . By (44) again, the convergence is
uniform on each bounded subset of $E$, so (44) passes to the
limit, and this limit $F$ satisfies the conclusion of
Theorem 2. This completes our proof.

\msi
{\bf References.}
\ssi
[A] P. Assouad, Plongements lipschitziens dans $\R^n$,
Bull. Soc. Math. France, 111(4), 429Ð448, 1983.
\ssi
[NN] Assaf Naor and Ofer Neiman,
Assouad's theorem with dimension independent of the snowflaking,
to be published, Revista Matematica Iberoamercana.

\bsi\bsi\bsi
Guy David, 
\pari
Univ Paris-Sud, Laboratoire de math\'{e}matiques UMR-8628, 
Orsay F-91405, France
\pari
and Institut Universitaire de France
\msi
Marie Snipes,
\pari
Department of Mathematics,
Hayes Hall,
Kenyon College,
Gambier, Ohio 43022, USA

\bye